\documentclass{elsart}

\usepackage{latexsym, amsmath, amssymb}
\usepackage{epic,eepic,epsfig}
\usepackage{amscd}
\usepackage{pstricks,pst-node,cite}

\numberwithin{equation}{section}

\ifx\blackandwhite\undefined
  \psset{linecolor=black}\else
  \psset{linecolor=black}
\fi

\psset{linewidth=.4pt,dash=3pt 3pt,doublesep=.05,dotsize=1pt 5}
\SpecialCoor

    \newcommand{\middlearrow}{\lput{:U}{\begin{pspicture}[0](0,0)(0,0)
\psline[arrows=->,arrowscale=1.5](2.2pt,0)(2.3pt,0)\end{pspicture}}}

\newcommand{\C}[1]{\mathcal #1}
\newcommand{\B}[1]{\mathbb #1}

\begin{document}
\begin{frontmatter}
\title{The chromatic numbers of double coverings of a graph}

\author{Dongseok Kim}
\address{Department of Mathematics \\ Kyungpook National University \\ Taegu, 702-201,
Korea} \ead{dongseok@knu.ac.kr}
\author{Jaeun Lee}
\address{Department of Mathematics\\ Yeungnam University\\ Kyongsan, 712-749, Korea}
\ead{julee@yu.ac.kr}
\thanks{The second author was supported in part by
Com$^2$MaC-KOSEF(R11-1999-054).}


\begin{abstract}
If we fix a spanning subgraph $H$ of a graph $G$, we can define a
chromatic number of $H$ with respect to $G$ and we show that it
coincides with the chromatic number of a double covering of $G$
with co-support $H$. We also find a few estimations for the
chromatic numbers of $H$ with respect to $G$.
\end{abstract}

\end{frontmatter}
\maketitle

\section{Introduction}

Let $G$ be a finite simple graph with vertex set $V(G)$ and
edge set $E(G)$. The cardinality of a set $X$ is denoted by $|X|$.
Throughout the paper, we assume all graphs are finite and simple.

The aim of this article is to find adequate formulae for the
chromatic numbers of covering graphs. The chromatic number
$\chi(G)$ of a graph $G$ is the smallest number of colors needed
to color the vertices of G so that no two adjacent vertices share
the same color. Since the exploratory paper by Dirac \cite{dirac},
the chromatic number has been in the center of graph theory
research. Its rich history can be found in several articles
\cite{HM, thmas}. The concept of covering graphs is relatively new
\cite{GTP, GTB}. Its precise definition can be given as follows.
For a graph $G$, we denote the set of all vertices adjacent to
$v\in V(G)$ by $N(v)$ and call it the \emph{neighborhood} of a
vertex $v$. A graph $\widetilde G$ is called a \emph{covering} of
$G$ with a projection $p:\widetilde G \to G$, if there is a
surjection $p:V(\widetilde G)\to V(G)$ such that $p|_{N(\tilde
v)}:N(\tilde v) \to N(v)$ is a bijection for any vertex $v\in
V(G)$ and $\tilde v \in p^{-1}(v)$. In particular, if $p$ is
two-to-one, then the projection $p:\widetilde G \to G$ is called a
\emph{double covering} of $G$. Some structures or properties of
graphs work nicely with covering graphs. The characteristic
polynomials of a covering graph $\widetilde G$ and its base graph
$G$ have a strong relation \cite{FKL:Character, KL:Character,
MS:Character}. The enumeration of non-isomorphic covering graphs
has been well studied \cite{Jones:Iclass, KL:Iclass}. Amit, Linial, and Matousek find
the asymptotic behavior of the chromatic numbers of $n$-fold coverings without
considering isomorphic types~\cite{ALM}. We will
relate the chromatic number and double covering graphs as follows.

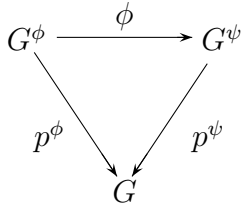
\begin{figure}
$$
\begin{pspicture}[.4](0,0.25)(3,2.75)
\qline(.6,2.5)(2.4,2.5)\psline[arrowscale=1.5]{->}(2.2,2.5)(2.4,2.5)
\qline(.3,2.3)(1.4,.65)\psline[arrowscale=1.5]{->}(1.2,.95)(1.4,.65)
\qline(2.6,2.15)(1.6,.65)\psline[arrowscale=1.5]{->}(1.8,.95)(1.6,.65)
\rput(.7,2.5){\rnode{a1}{$$}} \rput(.7,2.5){\rnode{a2}{$$}}
\rput(0.3,2.1){\rnode{a3}{$$}} \rput(2.7,2.1){\rnode{a4}{$$}}
\rput(0.2,2.5){\rnode{c1}{$G^\phi$}}
\rput(2.8,2.5){\rnode{c2}{$G^\psi$}}
\rput(1.5,0.45){\rnode{c3}{$G$}}
\rput[b](1.5,2.6){\rnode{c4}{$\phi$}}
\rput[tr](0.7,1.4){\rnode{c5}{$p^\phi$}}
\rput[tl](2.4,1.4){\rnode{c6}{$p^\psi$}}
\end{pspicture}
$$
\caption{Commuting diagram of isomorphic coverings.} \label{cd}
\end{figure}

A \emph{signed graph} is a pair $G_\phi = (G,~\phi)$ of a graph
$G$ and a function $\phi : E(G) \to {\B{Z}_2}$,
$\B{Z}_2=\{1,-1\}$. We call $G$ the \emph{underlying graph} of
$G_\phi$ and $\phi$ the \emph{signing} of $G$. A signing $\phi$ is
in fact a $\B{Z}_2$-voltage assignment of $G$, which was defined
by Gross and Tucker~\cite{GTP}.   It is known \cite{GTP, GTB} that
every double covering of a graph $G$ can be constructed as
follows: let $\phi$ be a signing of $G$. The double covering
$G^\phi$ of $G$ derived from $\phi$ has the following vertex set
$V(G^\phi)$ and edge set $E(G^\phi)$,
\begin{eqnarray*}
V(G^\phi)&=& \{ v_g | v\in V(G) ~\mathrm{and}~ g\in {\B{Z}_2}\} \\
E(G^\phi)&=& \{ (u_g, v_{\phi(u,v)g})| (u,v)\in E(G), g\in
{\B{Z}_2}\}.
\end{eqnarray*}
Two double coverings $p^\phi:G^\phi\to G $ and $p^\psi:G^\psi\to
G$ are \emph{isomorphic} if there exists a graph isomorphism $\phi
: G^\phi \to G^\psi$ such that the diagram in Figure \ref{cd}
commutes.

For a spanning subgraph $H$ of $G$, colorings $f$ and $g$ of $H$ are
\emph{compatible} in $G$ if for each edge $(u,v)\in E(G)-E(H)$,
$f(u)\not=g(v)$ and $f(v)\not=g(u)$. The smallest number of colors
such that $H$ has a pair of compatible colorings is called the
\emph{chromatic number of $H$ with respect to $G$} and denoted by
$\chi_{G}(H)$.  Since $(f|_{H},f|_{H})$ is a pair of compatible
colorings of $H$ for any spanning subgraph $H$ of $G$ and any
coloring $f$ of $G$, one can find $\chi_G(H) \leq \chi(G)$ for any
spanning subgraph $H$ of $G$. We remark that $\chi_{G}(G)=\chi(G)$
for any graph $G$, and that $\chi_{G}(\C{N}_{|V(G)|})=2$ if $G$ has
at least one edge, where $\C{N}_n$ is the null graph on $n$
vertices.

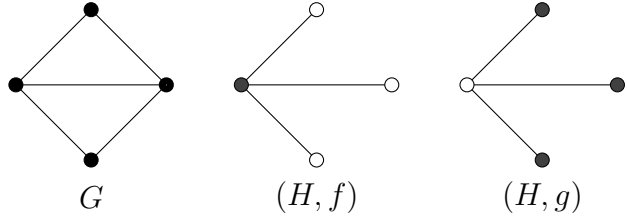
\begin{figure}
$$
\begin{pspicture}[.4](-1.8,-1.8)(1.5,1) \rput(.5,0){\rnode{a1}{$$}}
\rput(0,.5){\rnode{a2}{$$}} \rput(-.5,0){\rnode{a3}{$$}}
\rput(0,-.5){\rnode{a4}{$$}} \rput(1,0){\rnode{b1}{$$}}
\rput(.5,.5){\rnode{b2}{$$}} \rput(0,1){\rnode{b3}{$$}}
\rput(-.5,.5){\rnode{b4}{$$}} \rput(-1,0){\rnode{b5}{$$}}
\rput(-.5,-.5){\rnode{b6}{$$}} \rput(0,-1){\rnode{b7}{$$}}
\rput(.5,-.5){\rnode{b8}{$$}} \ncline{b1}{b3} \ncline{b1}{b5}
\ncline{b1}{b7} \ncline{b5}{b3} \ncline{b5}{b7}
\rput(0,-1.5){\rnode{c4}{$G$}} \pscircle[linewidth=2.5pt](1,0){.1}
\pscircle[fillstyle=solid,fillcolor=black](0,1){.1}
\pscircle[fillstyle=solid,fillcolor=black](-1,0){.1}
\pscircle[fillstyle=solid,fillcolor=black](0,-1){.1}
\end{pspicture}  \begin{pspicture}[.4](-1.5,-1.8)(4.5,1) \rput(.5,0){\rnode{a1}{$$}}
\rput(0,.5){\rnode{a2}{$$}} \rput(-.5,0){\rnode{a3}{$$}}
\rput(0,-.5){\rnode{a4}{$$}} \rput(1,0){\rnode{b1}{$$}}
\rput(.5,.5){\rnode{b2}{$$}} \rput(0,1){\rnode{b3}{$$}}
\rput(-.5,.5){\rnode{b4}{$$}} \rput(-1,0){\rnode{b5}{$$}}
\rput(-.5,-.5){\rnode{b6}{$$}} \rput(0,-1){\rnode{b7}{$$}}
\rput(.5,-.5){\rnode{b8}{$$}} \ncline{b5}{b7} \ncline{b5}{b3}
\ncline{b5}{b1} \rput(0,-1.5){\rnode{c4}{$(H,f)$}}
\pscircle[fillstyle=solid,fillcolor=white,linecolor=black](1,0){.1}
\pscircle[fillstyle=solid,fillcolor=white,linecolor=black](0,1){.1}
\pscircle[fillstyle=solid,fillcolor=darkgray,linecolor=black](-1,0){.1}
\pscircle[fillstyle=solid,fillcolor=white,linecolor=black](0,-1){.1}
\rput(3.5,0){\rnode{e1}{$$}} \rput(2,.5){\rnode{e2}{$$}}
\rput(2.5,0){\rnode{e3}{$$}} \rput(2,-.5){\rnode{e4}{$$}}
\rput(4,0){\rnode{f1}{$$}} \rput(3.5,.5){\rnode{f2}{$$}}
\rput(3,1){\rnode{f3}{$$}} \rput(2.5,.5){\rnode{f4}{$$}}
\rput(2,0){\rnode{f5}{$$}} \rput(2.5,-.5){\rnode{f6}{$$}}
\rput(3,-1){\rnode{f7}{$$}} \rput(3.5,-.5){\rnode{f8}{$$}}
\ncline{f5}{f7} \ncline{f5}{f3} \ncline{f5}{f1}
\rput(3,-1.5){\rnode{c4}{$(H,g)$}}
\pscircle[fillstyle=solid,fillcolor=darkgray,linecolor=black](4,0){.1}
\pscircle[fillstyle=solid,fillcolor=darkgray,linecolor=black](3,1){.1}
\pscircle[fillstyle=solid,fillcolor=white,linecolor=black](2,0){.1}
\pscircle[fillstyle=solid,fillcolor=darkgray,linecolor=black](3,-1){.1}
\end{pspicture}
$$
\caption{A spanning subgraph $H$ of $G$ with $\chi_{G}(H)=2$.}
\label{exam1}
\end{figure}

In Section~\ref{basic}, we recall some basic properties. We show that the
chromatic numbers of double coverings of a given graph can be computed from
the number $\chi_G(H)$ for any spanning subgraph $H$ of $G$. In
Section~\ref{estimation}, we will estimate the number $\chi_G(H)$.
We discuss a generalization to $n$-fold covering graphs in Section
\ref{fremark}.

\section{Basic properties}
\label{basic}

Let $\phi$ be a signing of $G$. We define the
\emph{support} of $\phi$ by the spanning subgraph of $G$ whose
edge set is $\phi^{-1}(-1)$, and denoted by $spt(\phi)$.
Similarly, we define the \emph{co-support} of $\phi$  by the
spanning subgraph of $G$ whose edge set is $\phi^{-1}(1)$, and
denoted by $cospt(\phi)$. Any spanning subgraph $H$ of $G$ can be
described as a co-support $cospt(\phi)$ of a signing $\phi$ of
$G$. Let $\phi_H$ be the signing of $G$ with $cospt(\phi_H)=H$.
Let $f$ and $g$ be compatible $\chi_G(H)$-colorings of $H$. We
define a function
$$h:V(G^\phi)\to \{1,2,\ldots,\chi_G(H)\}$$ by
$h(v_1)=f(v)$ and $h(v_{-1})=g(v)$ for each $v\in V(G)$. Then, by
the compatibility of $f$ and $g$, $h$ is a $\chi_G(H)$-coloring of
$G^\phi$. Hence, $\chi(G^\phi)\le \chi_G(H)$. Conversely, let $h$
be a $\chi(G^\phi)$-coloring of $G^\phi$. We define two
$\chi(G^\phi)$-colorings $f$ and $g$ of $H$ by $f(v)=h(v_1)$ and
$g(v)=h(v_{-1})$ for each $v\in V(G)$. Then $f$ and $g$ are
compatible because $h$ is a coloring of $G^\phi$. Hence,
$\chi_G(H)\le \chi(G^\phi)$. Now, we have the following theorem.

 \begin{thm}\label{doub}
 Let $H$ be a spanning subgraph of a graph $G$. Then
 $$\chi_G(H)=\chi(G^{\phi_H}),$$
 where $\phi_H$ is the signing of $G$ with
$cospt(\phi_H)=H$.
 \end{thm}

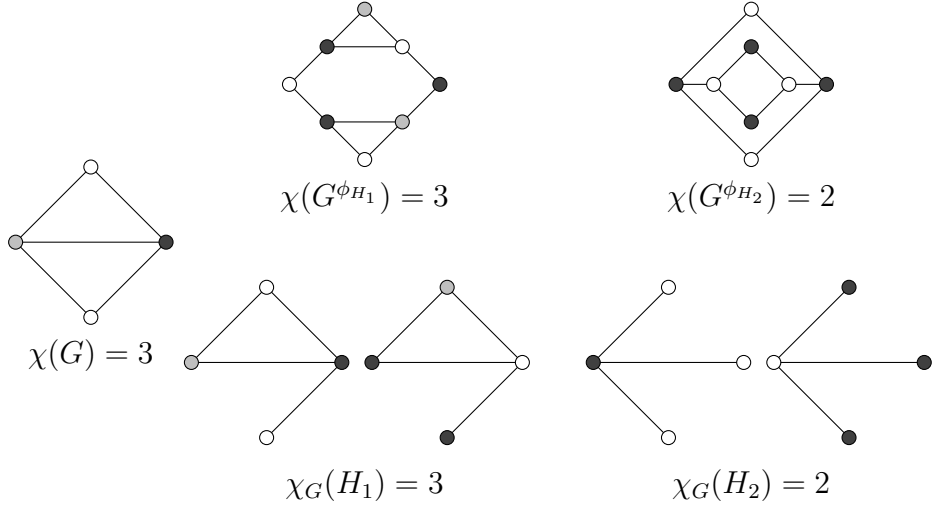
\begin{figure}
$$
\begin{matrix}
& \begin{pspicture}[.4](-1,-1.5)(1,1) \rput(.5,0){\rnode{a1}{$$}}
\rput(0,.5){\rnode{a2}{$$}} \rput(-.5,0){\rnode{a3}{$$}}
\rput(0,-.5){\rnode{a4}{$$}} \rput(1,0){\rnode{b1}{$$}}
\rput(.5,.5){\rnode{b2}{$$}} \rput(0,1){\rnode{b3}{$$}}
\rput(-.5,.5){\rnode{b4}{$$}} \rput(-1,0){\rnode{b5}{$$}}
\rput(-.5,-.5){\rnode{b6}{$$}} \rput(0,-1){\rnode{b7}{$$}}
\rput(.5,-.5){\rnode{b8}{$$}} \ncline{b1}{b3} \ncline{b1}{b7}
\ncline{b5}{b3} \ncline{b5}{b7} \ncline{b2}{b4} \ncline{b6}{b8}
\rput(0,-1.5){\rnode{c4}{$\chi(G^{\phi_{H_1}})=3$}}
\pscircle[fillstyle=solid,fillcolor=darkgray,linecolor=black](1,0){.1}
\pscircle[fillstyle=solid,fillcolor=lightgray,linecolor=black](0,1){.1}
\pscircle[fillstyle=solid,fillcolor=white,linecolor=black](-1,0){.1}
\pscircle[fillstyle=solid,fillcolor=white,linecolor=black](0,-1){.1}
\pscircle[fillstyle=solid,fillcolor=white,linecolor=black](.5,.5){.1}
\pscircle[fillstyle=solid,fillcolor=darkgray,linecolor=black](-.5,.5){.1}
\pscircle[fillstyle=solid,fillcolor=darkgray,linecolor=black](-.5,-.5){.1}
\pscircle[fillstyle=solid,fillcolor=lightgray,linecolor=black](.5,-.5){.1}
\end{pspicture} & \begin{pspicture}[.4](-1,-1.5)(1,1)
\rput(.5,0){\rnode{a1}{$$}} \rput(0,.5){\rnode{a2}{$$}}
\rput(-.5,0){\rnode{a3}{$$}} \rput(0,-.5){\rnode{a4}{$$}}
\rput(1,0){\rnode{b1}{$$}} \rput(.5,.5){\rnode{b2}{$$}}
\rput(0,1){\rnode{b3}{$$}} \rput(-.5,.5){\rnode{b4}{$$}}
\rput(-1,0){\rnode{b5}{$$}} \rput(-.5,-.5){\rnode{b6}{$$}}
\rput(0,-1){\rnode{b7}{$$}} \rput(.5,-.5){\rnode{b8}{$$}}
\ncline{b1}{a1} \ncline{b1}{b3} \ncline{b1}{b7} \ncline{b5}{b3}
\ncline{b5}{b7} \ncline{b5}{a3} \ncline{a1}{a2} \ncline{a2}{a3}
\ncline{a3}{a4} \ncline{a1}{a4}
\rput(0,-1.5){\rnode{c4}{$\chi(G^{\phi_{H_2}})=2$}}
\pscircle[fillstyle=solid,fillcolor=darkgray,linecolor=black](1,0){.1}
\pscircle[fillstyle=solid,fillcolor=white,linecolor=black](0,1){.1}
\pscircle[fillstyle=solid,fillcolor=darkgray,linecolor=black](-1,0){.1}
\pscircle[fillstyle=solid,fillcolor=white,linecolor=black](0,-1){.1}
\pscircle[fillstyle=solid,fillcolor=white,linecolor=black](.5,0){.1}
\pscircle[fillstyle=solid,fillcolor=darkgray,linecolor=black](0,.5){.1}
\pscircle[fillstyle=solid,fillcolor=white,linecolor=black](-.5,0){.1}
\pscircle[fillstyle=solid,fillcolor=darkgray,linecolor=black](0,-.5){.1}
\end{pspicture}  \\
\begin{pspicture}[.4](-1,-.6)(1,.6) \rput(.5,0){\rnode{a1}{$$}}
\rput(0,.5){\rnode{a2}{$$}} \rput(-.5,0){\rnode{a3}{$$}}
\rput(0,-.5){\rnode{a4}{$$}} \rput(1,0){\rnode{b1}{$$}}
\rput(.5,.5){\rnode{b2}{$$}} \rput(0,1){\rnode{b3}{$$}}
\rput(-.5,.5){\rnode{b4}{$$}} \rput(-1,0){\rnode{b5}{$$}}
\rput(-.5,-.5){\rnode{b6}{$$}} \rput(0,-1){\rnode{b7}{$$}}
\rput(.5,-.5){\rnode{b8}{$$}} \ncline{b1}{b3} \ncline{b1}{b5}
\ncline{b1}{b7} \ncline{b5}{b3} \ncline{b5}{b7}
\rput(0,-1.5){\rnode{c4}{$\chi(G)=3$}}
\pscircle[fillstyle=solid,fillcolor=darkgray,linecolor=black](1,0){.1}
\pscircle[fillstyle=solid,fillcolor=white,linecolor=black](0,1){.1}
\pscircle[fillstyle=solid,fillcolor=lightgray,linecolor=black](-1,0){.1}
\pscircle[fillstyle=solid,fillcolor=white,linecolor=black](0,-1){.1}
\end{pspicture}  & & \\
&\begin{pspicture}[.4](-1.2,-1.2)(1.2,1) \rput(.5,0){\rnode{a1}{$$}}
\rput(0,.5){\rnode{a2}{$$}} \rput(-.5,0){\rnode{a3}{$$}}
\rput(0,-.5){\rnode{a4}{$$}} \rput(1,0){\rnode{b1}{$$}}
\rput(.5,.5){\rnode{b2}{$$}} \rput(0,1){\rnode{b3}{$$}}
\rput(-.5,.5){\rnode{b4}{$$}} \rput(-1,0){\rnode{b5}{$$}}
\rput(-.5,-.5){\rnode{b6}{$$}} \rput(0,-1){\rnode{b7}{$$}}
\rput(.5,-.5){\rnode{b8}{$$}} \ncline{b1}{b3} \ncline{b1}{b7}
\ncline{b5}{b3} \ncline{b5}{b1}
\pscircle[fillstyle=solid,fillcolor=darkgray,linecolor=black](1,0){.1}
\pscircle[fillstyle=solid,fillcolor=white,linecolor=black](0,1){.1}
\pscircle[fillstyle=solid,fillcolor=lightgray,linecolor=black](-1,0){.1}
\pscircle[fillstyle=solid,fillcolor=white,linecolor=black](0,-1){.1}
\end{pspicture} \begin{pspicture}[.4](-1.2,-1.2)(1.4,1) \rput(.5,0){\rnode{a1}{$$}}
\rput(0,.5){\rnode{a2}{$$}} \rput(-.5,0){\rnode{a3}{$$}}
\rput(0,-.5){\rnode{a4}{$$}} \rput(1,0){\rnode{b1}{$$}}
\rput(.5,.5){\rnode{b2}{$$}} \rput(0,1){\rnode{b3}{$$}}
\rput(-.5,.5){\rnode{b4}{$$}} \rput(-1,0){\rnode{b5}{$$}}
\rput(-.5,-.5){\rnode{b6}{$$}} \rput(0,-1){\rnode{b7}{$$}}
\rput(.5,-.5){\rnode{b8}{$$}} \ncline{b1}{b3} \ncline{b1}{b7}
\ncline{b5}{b3} \ncline{b5}{b1}
\pscircle[fillstyle=solid,fillcolor=white,linecolor=black](1,0){.1}
\pscircle[fillstyle=solid,fillcolor=lightgray,linecolor=black](0,1){.1}
\pscircle[fillstyle=solid,fillcolor=darkgray,linecolor=black](-1,0){.1}
\pscircle[fillstyle=solid,fillcolor=darkgray,linecolor=black](0,-1){.1}
\end{pspicture} & \begin{pspicture}[.4](-1.4,-1.2)(1.2,1) \rput(.5,0){\rnode{a1}{$$}}
\rput(0,.5){\rnode{a2}{$$}} \rput(-.5,0){\rnode{a3}{$$}}
\rput(0,-.5){\rnode{a4}{$$}} \rput(1,0){\rnode{b1}{$$}}
\rput(.5,.5){\rnode{b2}{$$}} \rput(0,1){\rnode{b3}{$$}}
\rput(-.5,.5){\rnode{b4}{$$}} \rput(-1,0){\rnode{b5}{$$}}
\rput(-.5,-.5){\rnode{b6}{$$}} \rput(0,-1){\rnode{b7}{$$}}
\rput(.5,-.5){\rnode{b8}{$$}} \ncline{b5}{b7} \ncline{b5}{b3}
\ncline{b5}{b1}
\pscircle[fillstyle=solid,fillcolor=white,linecolor=black](1,0){.1}
\pscircle[fillstyle=solid,fillcolor=white,linecolor=black](0,1){.1}
\pscircle[fillstyle=solid,fillcolor=darkgray,linecolor=black](-1,0){.1}
\pscircle[fillstyle=solid,fillcolor=white,linecolor=black](0,-1){.1}
\end{pspicture} \begin{pspicture}[.4](-1.2,-1.2)(1.2,1) \rput(.5,0){\rnode{a1}{$$}}
\rput(0,.5){\rnode{a2}{$$}} \rput(-.5,0){\rnode{a3}{$$}}
\rput(0,-.5){\rnode{a4}{$$}} \rput(1,0){\rnode{b1}{$$}}
\rput(.5,.5){\rnode{b2}{$$}} \rput(0,1){\rnode{b3}{$$}}
\rput(-.5,.5){\rnode{b4}{$$}} \rput(-1,0){\rnode{b5}{$$}}
\rput(-.5,-.5){\rnode{b6}{$$}} \rput(0,-1){\rnode{b7}{$$}}
\rput(.5,-.5){\rnode{b8}{$$}} \ncline{b5}{b7} \ncline{b5}{b3}
\ncline{b5}{b1}
\pscircle[fillstyle=solid,fillcolor=darkgray,linecolor=black](1,0){.1}
\pscircle[fillstyle=solid,fillcolor=darkgray,linecolor=black](0,1){.1}
\pscircle[fillstyle=solid,fillcolor=white,linecolor=black](-1,0){.1}
\pscircle[fillstyle=solid,fillcolor=darkgray,linecolor=black](0,-1){.1}
\end{pspicture}\\
 & \chi_{G}(H_1)=3 &  \chi_{G}(H_2)=2
\end{matrix}
$$
\caption{$\chi_{G}(H)$ and the chromatic number of double coverings $G^{\phi_H}$ of a
graph $G$.} \label{exam2}
\end{figure}

It is not hard to see that the graph $G$ in Figure \ref{exam1} has
two non-isomorphic connected double coverings. We exhibit spanning
subgraphs $H_1, H_2$ of $G$ corresponding to two non-isomorphic
connected covering graphs of $G$ and their chromatic numbers with
respect to $G$ in Figure~\ref{exam2}.

For a subset $X \subset V(G)$ and for a spanning subgraph $H$ of
$G$, let $H_{X}$ denote a new spanning subgraph of $G$ defined as
follow: Two vertices in $X$ or in $V(G)-X$ are adjacent in $H_{X}$
if they are adjacent in $H$, while a vertex in $X$ and a vertex in
$V(G) - X $ are adjacent in $H_{X}$ if they are adjacent in $G$ but
not adjacent in $H$, i.e., they are adjacent in the complement
$\overline{H}(G)$ of $H$ in $G$. Two spanning subgraphs $H$ and $K$
of $G$ are \emph{Seidel switching equivalent} in $G$ if there exists
a subset $X \subset V(G)$ such that $H_{X}=K$. Clearly, the Seidel
switching equivalence is an equivalence relation on the set of
spanning subgraphs of $G$, and the equivalence class $[H]$ of a
spanning subgraph $H$ of $G$ is $\{H_{X}: X \subset V(G)\}$.

For a signing $\phi :E(G) \to {\B{Z}}_2$ and for any $X \subset
V(G)$, let $\phi_X$ be the signing obtained from $\phi$ by
reversing the sign of each edge having exactly one end point in
$X.$  If $\psi = \phi_X$ for some $X\subset V(G)$ then $\phi$ and
$\psi$ are said to be \emph{switching equivalent} \cite{CW}.

It is clear that for a subset $X \subset V(G)$ and for a spanning
subgraph $H$ of $G$, $H_X= cospt((\phi_H)_X)$. By a slight
modification of the proof of Corollary $4$~\cite{KHLS}, we obtain
the following theorem.

\begin{thm}\label{switchequi} Let $G$ be a graph.
Let $H, K$ be spanning subgraphs of $G$. Then the following
statements are equivalent.
\begin{description}
\item[{\rm (1)}] Two graphs $H$ and $K$ are Seidel switching
equivalent. \item[{\rm (2)}] Two signings $\phi_H$ and $\phi_K$
are switching equivalent. \item[{\rm (3)}] Two double coverings
$G^{\phi_H}$ and $G^{\phi_K}$ of $G$ are isomorphic as coverings.
\end{description}
\end{thm}

The following corollary follows easily from Theorem \ref{doub}
and \ref{switchequi}.

\begin{cor}\label{switchcon} Let $H$ and $K$ be two spanning subgraphs of a graph $G$. If
they are switching equivalent, then
$$\chi_G(H) =\chi_G(K).$$
\end{cor}

\begin{figure}
$$ \begin{pspicture}[.4](-1.8,-1.8)(1.5,1)  \rput(1,0){\rnode{b1}{$$}}
\rput(0,1){\rnode{b3}{$$}}
\rput(-1,0){\rnode{b5}{$$}}
\rput(0,-1){\rnode{b7}{$$}}
\ncline{b1}{b3} \ncline{b1}{b5}
\ncline{b1}{b7} \ncline{b5}{b3} \ncline{b5}{b7}
\rput(0,-1.5){\rnode{c4}{$G$}} \pscircle[linewidth=2.5pt](1,0){.1}
\pscircle[fillstyle=solid,fillcolor=black](0,1){.1}
\pscircle[fillstyle=solid,fillcolor=black](-1,0){.1}
\pscircle[fillstyle=solid,fillcolor=black](0,-1){.1}
\end{pspicture} \hskip 1cm \begin{pspicture}[.4](-1.8,-1.8)(1.5,1)  \rput(1,0){\rnode{b1}{$$}}
\rput(0,1){\rnode{b3}{$$}}
\rput(-1,0){\rnode{b5}{$$}}
\rput(0,-1){\rnode{b7}{$$}}
 \ncline{b1}{b5}
\ncline{b1}{b7} \ncline{b5}{b3} \ncline{b5}{b7}
\rput(0,-1.5){\rnode{c4}{$H$}} \pscircle[linewidth=2.5pt](1,0){.1}
\pscircle[fillstyle=solid,fillcolor=black](0,1){.1}
\pscircle[fillstyle=solid,fillcolor=black](-1,0){.1}
\pscircle[fillstyle=solid,fillcolor=black](0,-1){.1}
\end{pspicture} \hskip 1cm \begin{pspicture}[.4](-1.8,-1.8)(1.5,1)  \rput(1,0){\rnode{b1}{$$}}
\rput(0,1){\rnode{b3}{$$}}
\rput(-1,0){\rnode{b5}{$$}}
\rput(0,-1){\rnode{b7}{$$}}
\ncline{b1}{b3} \ncline{b1}{b5}
 \ncline{b5}{b3} \ncline{b5}{b7}
\rput(0,-1.5){\rnode{c4}{$K$}} \pscircle[linewidth=2.5pt](1,0){.1}
\pscircle[fillstyle=solid,fillcolor=black](0,1){.1}
\pscircle[fillstyle=solid,fillcolor=black](-1,0){.1}
\pscircle[fillstyle=solid,fillcolor=black](0,-1){.1}
\end{pspicture}
$$
\caption{Two non-switching
equivalent subgraphs $H, K$ in $G$ with $\chi_G(H) =\chi_G(K)$.} \label{nonswitching}
\end{figure}
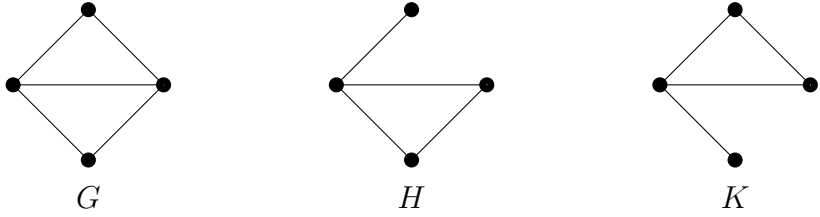

The converse of Corollary~\ref{switchcon} is not true in general. We have provided two non-switching
equivalent spanning subgraphs $H, K$ in $G$ with $\chi_G(H) =\chi_G(K)=3$ in Figure~\ref{nonswitching}.

For a coloring $f$ of $H$, let $\C{I}_f$ be the number of colors in
$\{1, 2, \ldots, \chi(H)\}$ such that the preimage $f^{-1}(i)$ is
independent in $\overline{H}(G)$ (and hence, also in $G$).

\begin{cor}\label{up-low-bds}
Let $G$ be a connected graph and let $H$ be a spanning subgraph of
$G$. Then
$$\max_{K \in [H]} \{\chi(K)\} \leq\chi_G(H) \leq \min_{K \in [H], f} \{\chi(G), 2\chi(K)-\C{I}_f \},$$
where  $f$ runs over all $\chi(K)$-colorings of $K$.
\end{cor}

\begin{pf}
It is clear that $\chi(H)\leq \chi_G(H)\leq \chi(G)$. Let $f$ be a
$\chi(H)$-coloring  of $H$ such that
$$\{ i | f^{-1}(i) ~\mathrm{is} ~\mathrm{independent}~\mathrm{in}~G \}
= \{\chi(H), \chi(H)-1, \ldots, \chi(H)-\C{I}_f +1\}.$$ We  define a
function $g:V(H)\to \{1,2,\ldots, 2\chi(H)-\C{I}_f\}$ as follows:
For a vertex $v$ in $V(H)$,
$$ g(v)=\left\{
\begin{array}{ll}
f(v) & \mbox{if $\chi(H)-\C{I}_f +1 \le f(v)\le \chi(H) $,} \\[1ex]
f(v)+\chi(H) & \mbox{otherwise.}\end{array}\right.$$ Then $g$ is a
coloring of $H$, and $f$ and $g$ are compatible. Now, the
corollary comes from Corollary~\ref{switchcon}.
\end{pf}

For a partition $\C{P}=\{V_1, V_2, \ldots, V_k\}$ of the vertex set
$V(G)$ of $G$, we define a new simple graph $G/\C{P}$ as follows:
the vertex set of $G/\C{P}$ is $\{V_1, V_2, \ldots, V_k\}$ and there
is an edge between two vertices $V_i$ and $V_j$ in $G/\C{P}$ if and
only if there exist two vertices $v_i\in V_i$ and $v_j\in V_j$ such
that $v_i$ and $v_j$ are adjacent in $G$ where $i\neq j$. We call $G/\C{P}$ the
\emph{quotient graph} associated with a partition $\C{P}$. For a
subset $S$ of $V(G)$, let $G[S]$ be the subgraph of $G$ whose vertex
set $S$ and whose edge set is the set of those edges of $G$ that
have both ends in $S$. We call $G[S]$ the \emph{subgraph induced by
$S$}.

\begin{cor}\label{h/p=bip}
Let $\C{P}=\{V_1,V_2, \ldots, V_k\}$ be a partition of the vertex
set of a connected graph $G$ and let $H=\cup_{i=1}^kG[V_i]$ be the
disjoint union of the induced subgraphs $G[V_i]$. If $G/\C{P}$ is
bipartite, then $\chi(G)=\chi_G(H).$
\end{cor}

\begin{pf} Let $X=\{[v_{i_1}], [v_{i_2}], \ldots, [v_{i_k}]\}$
be a part of the bipartition of the vertex set of the bipartite
graph $G/\C{P}$. Then $H_{\cup_{j=1}^{k} V_{i_j}}=G$. Then, the
corollary follows from Corollary ~\ref{switchcon}
or~\ref{up-low-bds}.
\end{pf}

The following theorem finds a necessary and sufficient condition for the bipartiteness of
covering graphs.

\begin{thm} [\cite{AKLS:BIPARTITE}]
Let $G$ be a non-bipartite graph with a generating voltage assignment
$\nu$ in $\mathcal{A}$ which derives the covering graph $\tilde{G}$.
Then $\tilde{G}$ is bipartite if and only if there
exists a subgroup $\mathcal{A}_e$ of index two in $\mathcal{A}$ such that for every cycle $C, \nu(C) \in \mathcal{A}_e$ if
and only if the length of $C$ is even. \label{thm45}
\end{thm}

It is obvious that $\chi_G(H)=1$ if and only if $G$ is a null graph.
In Theorem~\ref{bipar2}, we find a necessary and sufficient condition of $\chi_G(H)=2$.

\begin{thm}
Let $G$ be a connected graph having at least one edge and let $H$ be
a spanning subgraph of $G$. Then $\chi_G(H)=2$ if and only if either
$G$ is bipartite or $H\in [\C{N}_{|V(G)|}]$, where $\C{N}_{|V(G)|}$
is the null graph on $|V(G)|$ vertices. \label{bipar2}
\end{thm}

\begin{pf}
Let $G$ be a bipartite graph and $H$ be a spanning subgraph of
$G$. Then there exists a graph $K$ in the switching class $[H]$ of
$H$ in $G$ such that $K$ has at least one edge. By
Corollary~\ref{up-low-bds}, $\chi_G(H)=2$. We recall that $G$
itself is a spanning subgraph of $G$ and $\chi_G(G)=\chi(G)$.
Therefore, if a graph $G$ has at least one edge, then $G$ is
bipartite if and only if $\chi_G(H)=2$ for any spanning subgraph
$H$ of $G$.

If $G$ is not a bipartite graph, there exist a
signing $\phi$ of $G$ such that $G^\phi$ is bipartite~\cite{GTP}. It follows the
connectedness of $G$ that there exists a subset $Y$ of $V(G)$ such
that $cospt(\phi_Y)$ is connected. Since $G^\phi$ and $G^{\phi_Y}$
are isomorphic, by Theorem~\ref{switchequi}, $G^{\phi_Y}$ is
bipartite. We note that $cospt(\phi_Y)$ is isomorphic to a
subgraph of $G^{\phi_Y}$ and hence it is bipartite. Let $e$ be an
edge of $G$ such that one end is in $Y$ and the other is in
$V(G)-Y$. If $\phi_Y(e)=-1$, then there exists an even cycle which
contains the edge $e$ as the only edge whose value under $\phi_Y$
is $-1$. It follows from Theorem~\ref{thm45} that
$G^{\phi_Y}$ is not bipartite. This is a contradiction. It implies
that for an edge $e$ of $G$, $\phi_Y(e)=1$ if and only if one end
of $e$ is in $Y$ and the other is in $V(G)-Y$. Let $X$ be a part
of the bipartition of $cospt(\phi_Y)$, i.e., every edge $e$ in
$cospt(\phi_Y)$ has one end in $X$ and the other end in $V(G)-X$.
Then $cospt((\phi_{Y})_{_X})=\C{N}_{|V(G)|}$. Notice that
$cospt((\phi_{Y})_{_X})=cospt(\phi_Z)$, where $Z=(Y-X)\cup (X-Y)$.
Hence $H$ is switching equivalent to $\C{N}_{|V(G)|}$. It completes the proof of theorem.
\end{pf}

\section{Computations of $\chi_G(H)$} \label{estimation}

In this section, we aim to estimate the number $\chi_G(H)$ for any
spanning subgraph $H$ of $G$. Let $H$ be a spanning subgraph of
$G$, and let $F_1, F_2, \ldots, F_k$ be the components of the
complement $\overline{H}(G)$ of $H$ in $G$. Then
$\C{P}_{\bar{H}}=\{V(F_1)$, $V(F_2)$, $\ldots$ , $V(F_k)\}$ is a partition
of the vertex set $V(H)=V(G)$. Now, we consider a
$\chi(H/\C{P}_{\bar{H}})$-coloring $c$ of the quotient graph
$H/\C{P}_{\bar{H}}$. Then $c$ induces
a partition
$\C{P}_c=\{c^{-1}(1)$, $\ldots$, $c^{-1}$
$(\chi(H/\C{P}_{\bar{H}}))\}$ of the vertex set
$H/\C{P}_{\bar{H}}$. By composing the quotient map $: G\rightarrow H/\C{P}_{\bar{H}}$ and $c:H/\C{P}_{\bar{H}} \rightarrow \{1, 2, \ldots,\chi(H/\C{P}_{\bar{H}})\}$, we obtain a partition of $H$ and by slightly abusing notation we denoted it identically $\C{P}_c$. One can notice that each vertex of
$H/\C{P}_c$ can be considered as a union of the vertex sets
$V(F_1), \ldots, V(F_k)$. For each $i=1, \ldots,
\chi(H/\C{P}_{\bar{H}})$, let $H_c(i)=H[c^{-1}(i)]$, where we
consider $c^{-1}(i)$ as a subset of $V(H)=V(G)$. A coloring $f$ of
$H$ \emph{respects the coloring $c$} of $H/\C{P}_{\bar{H}}$ if
$|f(H_c(i))|=\chi(H_c(i))$ and $f(H_c(i))\cap f(H_c(j))=\emptyset$
for any $1\le i\not=j\le\chi(H/\C{P}_{\bar{H}})$. For a coloring
$f$ which respects $c$, let $\C{I}_f(i)$ be the number of colors
in $\{i_1, i_2, \ldots, i_{\chi(H_c(i))}\}$ such that the vertex
set $f^{-1}(i_k)$ is independent in $\overline{H}(G)$ and let
$\C{D}_f(i)=\chi(H_c(i))-\C{I}_f(i)$ for each $i=1, \ldots,
\chi(H/\C{P}_{\bar{H}})$.  Let
$$\Delta_S=  \max\left\{0,\, 2\,\max\{s\,|\, s\in S\}-\sum_{s\in S} s \right\}$$
for any subset $S$ of natural numbers.

\begin{thm}\label{complem}
Let $G$ be a connected graph  and let $H$ be a spanning subgraph of
$G$. Then
$$\begin{array}{lcl}
\chi_G(H) & \le & \displaystyle \min_c\left\{ \displaystyle
\sum_{i=1}^{\chi(H/\C{P}_{\bar{H}})}\!\!\! \chi(H_c(i)) +
\Delta_{\{\chi(H_c(i))\,|\,i=1,2, \ldots,
\chi(H/\C{P}_{\bar{H}})\}},\right.\\[3ex]
& & \hspace{1.5cm} \left.\displaystyle
\sum_{i=1}^{\chi(H/\C{P}_{\bar{H}})}\!\!\! \chi(H_c(i))+
\min_{f}\left\{\Delta_{\{\C{D}_f(i)\,|\,i=1,2, \ldots,
\chi(H/\C{P}_{\bar{H}})\}} \right\}\right\},\end{array}$$ where $c$
runs over all ${\chi(H/\C{P}_{\bar{H}})}$-colorings of
$H/\C{P}_{\bar{H}}$ and $f$ runs over all colorings of $H$ which
respect $c$.
\end{thm}

\begin{pf}
Let $c$ be a $\chi(H/\C{P}_{\bar{H}})$-coloring of
$H/\C{P}_{\bar{H}}$ and let $f$ be a coloring of $H$ which respects
$c$.

 First, we want to show that
$$\chi_G(H)\le\sum_{i=1}^{\chi(H/\C{P}_{\bar{H}})} \chi(H_c(i)) +
\Delta_{\{\chi(H_c(i))\,|\,i=1,2, \ldots,
\chi(H/\C{P}_{\bar{H}})\}}.$$

Without loss of generality, we may assume that
$\chi(H_c(1))\ge \chi(H_c(2))\ge $ $\ldots $ $\ge$
$\chi(H_c(\chi(H/\C{P}_{\bar{H}})))$.
Let the
image $$f(V(H_c(i)))=\{\sum_{j=1}^{i-1}\chi(H_c(j))+1, \ldots,
\sum_{j=1}^i\chi(H_c(j))\}$$ and
$$\ell=\Delta_{\{\chi(H_c(i))\})\,|\,i=1,2, \ldots,
\chi(H/\C{P}_{\bar{H}})\}}.$$ Then
$$\ell=\max\{\chi(H_c(1))-\sum_{i=2}^{\chi(H/\C{P}_{\bar{H}})}\chi(H_c(i)),\,0\}.$$
We define $g:V(H) \to \{1,2,\ldots, n+\ell\}$ by
$g(v)=f(v)-\chi(H_c(1))$, where
$$n=\sum_{i=1}^{\chi(H/\C{P}_{\bar{H}})}\chi(H_c(i))$$ and the
arithmetic is done by modulo $n+\ell$. Then $g$ is a coloring of
$H$. Since $f(V(H_c(i)))\cap g(V(H_c(i)))=\emptyset$ and each edge
in $E(G)-E(H)=E(\overline{H}(G))$ connects two vertices in $H_c(i)$
for some $i=1,2,\ldots, \chi(H/\C{P}_{\bar{H}})$, we can see that
$f$ and $g$ are compatible. Hence, $$\chi_G(H)\le
\sum_{i=1}^{\chi(H/\C{P}_{\bar{H}})}\chi(H_c(i))+\ell.$$

Next, we want to show that
$$\chi_G(H)\le\sum_{i=1}^{\chi(H/\C{P}_{\bar{H}})} \chi(H_c(i)) +
\Delta_{\{\C{D}_f(i)\,|\,i=1,2, \ldots,
\chi(H/\C{P}_{\bar{H}})\}}.$$ In general, we may assume that
$\C{D}_f(1)\ge \C{D}_f(2)\ge \cdots \ge
\C{D}_f(\chi(H/\C{P}_{\bar{H}}))$. Now, we aim to define another
coloring $g$ of $H$ such that $f$ and $g$ are compatible. To do
this, first for the vertices $v$ of $H$ such that the set
$f^{-1}(f(v))$ is independent in $\overline{H}(G)$, we define
$g(v)=f(v)$. Next, by using a method similar to the first case, we
can extend the function $g$ to whole graph $H$ so that $f$ and $g$
are  compatible colorings of  $H$.

Finally, by taking the minimum value among all
$\chi(H/\C{P}_{\bar{H}})$-coloring $c$ of $H/\C{P}_{\bar{H}}$ and
all coloring $f$ of $H$ which respect $c$, we have the theorem.
\end{pf}

The following example shows the upper bound in Theorem~\ref{complem}
is sharp.
\begin{exmp}
Let $m,n$ be integers with $2\le m \le n$. Let $K_{m-1}$ be the
complete graph on $m-1$ vertices $v_1,\ldots, v_{m-1}$. Let $H_m$ be
a spanning subgraph of $K_{n}$ obtained by adding $n-m+1$ isolated
vertices $v_{m},\ldots, v_n$ to $K_{m-1}$. Then $\chi_{K_n}(H_m)=m$.
\label{example32}
\end{exmp}
\begin{pf}
To show $m\le \chi_{K_n}(H_m)$, we set $X=V(K_{m-1})$. Then
$\chi((H_m)_X)=m$ and hence $m\le \chi_{K_n}(H_m)$ by
Corollary~\ref{up-low-bds}. We can show that $\chi_{K_n}(H_m)\le m$
by using two methods which are contained in the proof of
Theorem~\ref{complem}. For the first method, we replace $H_m$ by
$(H_m)_X$. We observe that
$\overline{(H_m)_X}(K_n)=K_{n-m+1}\cup\{v_1,\ldots, v_{m-1}\}$ and
$(H_m)_X/\C{P}_{\bar{(H_m)}_X}=K_{m}$. Let $c$ be a $(m)$-coloring
of $(H_m)_X/\C{P}_{\bar{H}_X}$ such that $V((H_m)_c(i))=\{v_i\}$ for
each $i=1,2,\ldots, m-1$ and  $V((H_m)_c(m+1))=\{v_{m}, \ldots,
v_n\}$. We note that $\chi((H_m)_c(i))=1$ for each $i=1,2,\ldots,
m$. Since $\Delta_{\{1,1,\ldots,1\}}=0$, by Theorem~\ref{complem},
we have $\chi_{K_n}(H_m)=\chi_{K_n}((H_m)_X)\le m$. For the second
method, let $c$ be the trivial coloring of
$H_m/\C{P}_{\overline{H_m}}=K_1$ and let $f$ be a $(m-1)$-coloring
of $H_m$ such that $f(v_i)=i$ for each $i=1,2,\ldots, m-1$ and
$f(v_{m})=f(v_{m+1})=\cdots =f(v_n)=1$. Then $f$ respects $c$ and
$$\C{D}_f(1)=\chi((H_m)_c(1))-\C{D}_f(1)=\chi(H_m)-\C{D}_f(1)=(m-1)-(m-2)=1.$$
Since $\Delta_{\{1\}}=2-1=1$, by  Theorem~\ref{complem}, we have
$\chi_{K_n}(H_m)\le m$.
\end{pf}

Example \ref{example32} can be generalized to the following
corollary.

\begin{cor}\label{kpart}
Let $H$ be a complete $m$-partite graph which is a spanning subgraph
of $K_n$. Then $\chi_{K_n}(H)=m$.
\end{cor}

\begin{pf}
We observe that the complement $\overline{H}(K_n)$ of $H$ is also a
spanning subgraph of $K_n$ having at least $k$ components of which
each vertex set is a subset of a part of $H$. It is not hard to show
that $H/\C{P}_{\bar{H}}$ is also a complete $m$-partite graph and
$\chi(H_c(i))=1$ for each $i=1,\ldots, m$. By Theorem~\ref{complem},
$\chi_G(H)\le m$. Since $\chi(H)=m$, by Corollary~\ref{up-low-bds},
it completes the proof.
\end{pf}

If each component of a spanning subgraph $H$ of a graph $G$ is a
vertex induced subgraph, we can have an upper bound of the chromatic
number of $H$ induced by $G$ which is simpler than that in
Theorem~\ref{complem}.

\begin{thm}\label{part}
Let $\C{P}=\{V_1,V_2, \ldots, V_k\}$ be a partition of the  vertex
set of a connected graph $G$. Let $H$ be a disjoint union of the
induced subgraphs $G[V_1]$, $G[V_2]$,$ \ldots$, $G[V_n]$. Then we
have
$$\max_{V_i, V_j}\{ \chi(G[V_i\cup V_j])\}\le \chi_G(H)\le
\max_{V_i, V_j}\{ \chi(G[V_i])+\chi(G[V_j])\},$$ where  $V_i$ and
$V_j$ runs over all pairs of adjacent vertices in
 $G/\C{P}$.
\end{thm}

\begin{pf}
Let $V_i$ and $V_j$ be two adjacent vertices in $G/\C{P}$. Then
$G[V_i\cup V_j]$ is a subgraph of $H_{V_i}$. By
Corollary~\ref{up-low-bds}, $\chi(G[V_i\cup V_j])\le \chi(H_{V_i})
\le \chi_G(H)$ and hence
$$\max\{\chi(G[V_i\cup V_j])\,|\, \mbox{\rm $V_i$ is
adjacent to $V_j$ in  $G/\C{P}$}\}\le \chi_G(H).$$ For the second
inequality, let
$$M= \max\{ \chi(G[V_i])+\chi(G[V_j])\,|\,
V_i ~\mathrm{is}~\mathrm{adjacent}~\mathrm{to}~V_j~\mathrm{in}
~G/\C{P}\}.$$

By the definition of $M$, there exist $s$, $t$ and $M$-coloring $f:V(H)\to$ $\{1$, $2$, $\ldots$, $M\}$ of $H$ such
that
$$ \chi(G[V_s])+\chi(G[V_t])=M,$$
and $$f(G[V_i])=\{1,2,\ldots,\chi(G[V_i])\}$$ for each $i=1$, $2$,
$\ldots$, $k$. We note that $f$ may not be surjective. We define
another $M$-coloring $g$ of $H$  such that $g(G[V_i])=\{M, M-1,
\ldots, M-\chi(G[V_i])+1\}$. Now, we aim to show that $f$ and $g$
are compatible.
 Let $uv$
be an edge of $G$ which is not in $E(H)$.  Now, by the construction
of $G/\C{P}$, Then there exist $i$ and $j$ such that $u\in V_i$ and
$v\in V_j$. By the definition of $G/\C{P}$, $V_i$ is adjacent to
$V_j$ in $G/\C{P}$. If $f(V_i)\cap g(V_j)\not=\emptyset$, then, by
the construction of $f$ and $g$, $M<\chi(G[V_i])+\chi(G[V_j])$. This
contradicts the hypothesis of $M$. Thus, $f(V_i)\cap
g(V_j)=\emptyset$. Similarly, we can see that $g(V_i)\cap
f(V_j)=\emptyset$. Therefore, $f(u)\not=g(v)$ and $g(u)\not=f(v)$,
i.e., $f$ and $g$ are compatible. It completes the proof.
\end{pf}

By Corollary~\ref{h/p=bip} and Theorem~\ref{part}, we have the
following corollaries.

\begin{cor}\label{h/p=bip'}
Let $\C{P}=\{V_1,V_2, \ldots, V_k\}$ be a partition of the  vertex
set of a connected graph $G$.
 If $G/\C{P}$ is bipartite, then
 $$\max_{V_i, V_j}\{ \chi(G[V_i\cup V_j])\}\le \chi(G)\le
\max_{V_i, V_j}\{ \chi(G[V_i])+\chi(G[V_j])\},$$ where  $V_i$ and
$V_j$ runs over all pairs of adjacent vertices in
 $G/\C{P}$.
\end{cor}

\begin{cor}\label{kn1}
Let $H$ be a spanning subgraph of a connected graph $G$ such that $H$ has $k$
components $H_1, H_2, \ldots, H_k$ with $\chi(H_i)\ge
\chi(H_{i+1})$ for each $i=1,2, \ldots, k-1$.
\begin{enumerate}
\item[{\rm (1)}]If the complement $\overline{H}(G)$ of $H$ in $G$
is the complete $k$ partite graph, then we have
$\chi_G(H)=\chi(H_1)+\chi(H_2)$.
\item[{\rm (2)}]If each component of $H$ is the complete graph,
$i.e$, $H_i=K_{\ell_i}$ for each $i=1,2, \ldots, k$. Then we have
$\chi_G(H)\le \ell_1+\ell_2$. In particular, if $G$ is the complete
graph $K_n$, then we have $\chi_G(H)=\ell_1+\ell_2$.
\end{enumerate}
\end{cor}

\begin{pf}
We observe that, in any case,  each component $H_i$ of $H$ is an
induced subgraph $G[V(H_i)]$ for each $i=1,2,\ldots,k$, and
$\C{P}=\{V(H_i)\,|\, i=1,2, \ldots, k\}$ forms a partition of $G$.

(1) If the complement $\overline{H}(G)$ is the
  complete $k$ partite graph, then $H/\C{P}$ is the complete graph
  $K_k$, i.e., each pair of vertices in $H/\C{P}$ is adjacent in
  $H/\C{P}$. Since $\chi(G[V(H_1)\cup V(H_2)])=\chi(H_1)+\chi(H_2)$
and $\max\{ \chi(G[V(H_i)])+\chi(G[V(H_j)])\,|\, 1\le i\not=j \le k
\}= \chi(H_1)+\chi(H_2)$,  by Theorem~\ref{part}, we have
  $\chi_G(H)=\chi(H_1)+\chi(H_2)$.

  (2) If $H_i=K_{\ell_i}$ for each $i=1,2,
\ldots, k$, then, By Theorem~\ref{part}, we have  $\chi_G(H)\le
\ell_1+\ell_2$. If $G$ is the complete graph $K_n$, then
$\overline{H}(G)$ of $H$ in $G$ is the
  complete $k$ partite graph, by (1), we have
  $\chi_G(H)=\ell_1+\ell_2$ \cite{AKLS:BIPARTITE}.
\end{pf}

\section{Further remarks}
\label{fremark}

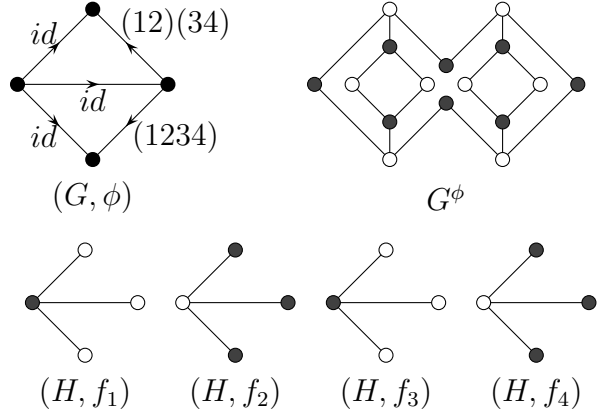
\begin{figure}
\begin{align*}
\begin{pspicture}[.4](-1.8,-1.8)(1.5,1) \rput(.5,0){\rnode{a1}{$$}}
\rput(0,.5){\rnode{a2}{$$}} \rput(-.5,0){\rnode{a3}{$$}}
\rput(0,-.5){\rnode{a4}{$$}} \rput(1,0){\rnode{b1}{$$}}
\rput(.5,.5){\rnode{b2}{$$}} \rput(0,1){\rnode{b3}{$$}}
\rput(-.5,.5){\rnode{b4}{$$}} \rput(-1,0){\rnode{b5}{$$}}
\rput(-.5,-.5){\rnode{b6}{$$}} \rput(0,-1){\rnode{b7}{$$}}
\rput(.5,-.5){\rnode{b8}{$$}} \ncline{b1}{b3}\middlearrow
\ncline{b5}{b1}\middlearrow \ncline{b1}{b7}\middlearrow
\ncline{b5}{b3}\middlearrow \ncline{b5}{b7}\middlearrow
\rput(0,-1.5){\rnode{c1}{$(G, \phi)$}}
\rput(0,-.2){\rnode{c2}{$id$}} \rput(-.65,.69){\rnode{c3}{$id$}}
\rput(-.65,-.69){\rnode{c4}{$id$}}
\rput(1.1,.8){\rnode{c5}{$(12)(34)$}}
\rput(1.1,-.7){\rnode{c6}{$(1234)$}}
\pscircle[linewidth=2.7pt](1,0){.1}
\pscircle[linewidth=2.7pt](0,1){.1}
\pscircle[linewidth=2.7pt](-1,0){.1}
\pscircle[linewidth=2.7pt](0,-1){.1}
\end{pspicture}
\begin{pspicture}[.4](-3.2,-1.8)(2.2,1)
\psline(-.75,1)(0,.25)(.75,1)(1.75,0)(.75,-1)(0,-.25)(-.75,-1)(-1.75,0)(-.75,1)(-.75,.5)(-1.25,0)(-.75,-.5)(-.75,-1)
\psline(-.75,.5)(-.25,0)(-.75,-.5)
\psline(.75,1)(.75,.5)(1.25,0)(.75,-.5)(.75,-1)
\psline(.75,.5)(.25,0)(.75,-.5)
\pscircle[fillstyle=solid,fillcolor=darkgray,linecolor=black](1.75,0){.1}
\pscircle[fillstyle=solid,fillcolor=darkgray,linecolor=black](.75,.5){.1}
\pscircle[fillstyle=solid,fillcolor=darkgray,linecolor=black](.75,-.5){.1}
\pscircle[fillstyle=solid,fillcolor=darkgray,linecolor=black](0,.25){.1}
\pscircle[fillstyle=solid,fillcolor=darkgray,linecolor=black](0,-.25){.1}
\pscircle[fillstyle=solid,fillcolor=darkgray,linecolor=black](-.75,.5){.1}
\pscircle[fillstyle=solid,fillcolor=darkgray,linecolor=black](-.75,-.5){.1}
\pscircle[fillstyle=solid,fillcolor=darkgray,linecolor=black](-1.75,0){.1}
\pscircle[fillstyle=solid,fillcolor=white,linecolor=black](-.75,1){.1}
\pscircle[fillstyle=solid,fillcolor=white,linecolor=black](-.75,-1){.1}
\pscircle[fillstyle=solid,fillcolor=white,linecolor=black](-1.25,0){.1}
\pscircle[fillstyle=solid,fillcolor=white,linecolor=black](-.25,0){.1}
\pscircle[fillstyle=solid,fillcolor=white,linecolor=black](.25,0){.1}
\pscircle[fillstyle=solid,fillcolor=white,linecolor=black](1.25,0){.1}
\pscircle[fillstyle=solid,fillcolor=white,linecolor=black](.75,1){.1}
\pscircle[fillstyle=solid,fillcolor=white,linecolor=black](.75,-1){.1}
\rput(0,-1.5){\rnode{c1}{$G^{\phi}$}}
\end{pspicture}
\\
\begin{pspicture}[.4](-1,-1.8)(1,1) \rput(.7,0){\rnode{a1}{$$}}
\rput(0,.7){\rnode{a2}{$$}} \rput(-.7,0){\rnode{a3}{$$}}
\rput(0,-.7){\rnode{a4}{$$}} \ncline{a3}{a2} \ncline{a3}{a1}
\ncline{a3}{a4} \rput(0,-1.2){\rnode{c4}{$(H,f_1)$}}
\pscircle[fillstyle=solid,fillcolor=white,linecolor=black](.7,0){.1}
\pscircle[fillstyle=solid,fillcolor=white,linecolor=black](0,.7){.1}
\pscircle[fillstyle=solid,fillcolor=darkgray,linecolor=black](-.7,0){.1}
\pscircle[fillstyle=solid,fillcolor=white,linecolor=black](0,-.7){.1}
\end{pspicture}
\begin{pspicture}[.4](-1,-1.8)(1,1) \rput(.7,0){\rnode{a1}{$$}}
\rput(0,.7){\rnode{a2}{$$}} \rput(-.7,0){\rnode{a3}{$$}}
\rput(0,-.7){\rnode{a4}{$$}} \ncline{a3}{a2} \ncline{a3}{a1}
\ncline{a3}{a4} \rput(0,-1.2){\rnode{c4}{$(H,f_2)$}}
\pscircle[fillstyle=solid,fillcolor=darkgray,linecolor=black](.7,0){.1}
\pscircle[fillstyle=solid,fillcolor=darkgray,linecolor=black](0,.7){.1}
\pscircle[fillstyle=solid,fillcolor=white,linecolor=black](-.7,0){.1}
\pscircle[fillstyle=solid,fillcolor=darkgray,linecolor=black](0,-.7){.1}
\end{pspicture}
\begin{pspicture}[.4](-1,-1.8)(1,1) \rput(.7,0){\rnode{a1}{$$}}
\rput(0,.7){\rnode{a2}{$$}} \rput(-.7,0){\rnode{a3}{$$}}
\rput(0,-.7){\rnode{a4}{$$}} \ncline{a3}{a2} \ncline{a3}{a1}
\ncline{a3}{a4} \rput(0,-1.2){\rnode{c4}{$(H,f_3)$}}
\pscircle[fillstyle=solid,fillcolor=white,linecolor=black](.7,0){.1}
\pscircle[fillstyle=solid,fillcolor=white,linecolor=black](0,.7){.1}
\pscircle[fillstyle=solid,fillcolor=darkgray,linecolor=black](-.7,0){.1}
\pscircle[fillstyle=solid,fillcolor=white,linecolor=black](0,-.7){.1}
\end{pspicture}
\begin{pspicture}[.4](-1,-1.8)(1,1) \rput(.7,0){\rnode{a1}{$$}}
\rput(0,.7){\rnode{a2}{$$}} \rput(-.7,0){\rnode{a3}{$$}}
\rput(0,-.7){\rnode{a4}{$$}} \ncline{a3}{a2} \ncline{a3}{a1}
\ncline{a3}{a4} \rput(0,-1.2){\rnode{c4}{$(H,f_4)$}}
\pscircle[fillstyle=solid,fillcolor=darkgray,linecolor=black](.7,0){.1}
\pscircle[fillstyle=solid,fillcolor=darkgray,linecolor=black](0,.7){.1}
\pscircle[fillstyle=solid,fillcolor=white,linecolor=black](-.7,0){.1}
\pscircle[fillstyle=solid,fillcolor=darkgray,linecolor=black](0,-.7){.1}
\end{pspicture}
\end{align*}
\caption{A permutation voltage assignment $\phi$ of $G$, compatible
colorings of $H = cospt(\phi)$ with $\chi_{G}(H)=2$ and its
corresponding $4$-fold covering graph $G^{\phi}$.} \label{exam3}
\end{figure}

\subsection{Existence of a spanning subgraph $H$ of $G$ with $\chi_G(H) = m$  for any $m$ with $2\le m \le \chi(G)$}

For $n\ge 2$, and for any spanning subgraph $H$ of a complete
graph $K_n$, we have $2\le \chi_{K_n}(H)\le n$. For converse, we
showed that for any integer $m$, between $2$ and $n$, there exists
a spanning subgraph $H_m$ of $K_n$ such that $\chi_{K_n}(H_m)=m$
in Example \ref{example32}. One can ask this can be extended to an
arbitrary connected graph. Let $G$ be connected graph. For any $m$ with $2\le m \le \chi(G)$, let $H$ be
an $m$-critical subgraph of $G$, that is, $\chi(H)=m$ and for any  proper subgraph $S$ of $H$ $\chi(S) <m$. Let $\tilde{H}$ be
a  spanning subgraph of $G$ obtained by adding $|\overline{H}(G)|$ isolated
vertices to $H$. By Theorem \ref{part}, $\chi_G(\tilde{H})=m$.

\subsection{n-fold covering graphs}

For $n$-fold covering graphs, let $S_n$ denote a symmetric group on
$n$ elements $\{1,2,\ldots, n\}$. Every edge of a graph $G$ gives
rise to a pair of oppositely directed edges. We denote the set of
directed edges of $G$ by $D(G)$. By $e^{-1}$ we mean the reverse
edge to an edge $e$. Each directed edge $e$ has an initial vertex
$i_e$ and a terminal vertex $t_e$. A {\it permutation voltage
assignment} $\phi$ on a graph $G$ is a map $\phi :D(G) \rightarrow
S_n$ with the property that $\phi(e^{-1})=\phi(e)^{-1}$ for each $e
\in D(G)$. The {\it permutation derived graph} $G^{\phi}$ is defined
as follows: $V(G^{\phi})=V(G) \times \{1, \dots ,n\}$, and for each
edge $e \in D(G)$ and $j \in \{1, \dots ,n\}$ let there be an edge
$(e,j)$ in $D(G^{\phi})$ with $i_{(e,j)} = (i_e,j)$ and
$t_{(e,j)}=(t_e,\phi(e)j)$. The natural projection $p_{\phi}
:G^{\phi} \rightarrow G$ is a covering. In \cite{GTP,GTB}, Gross and
Tucker showed that
 every $n$-fold covering ${\widetilde G}$ of a graph
$G$ can be derived from a voltage assignment.

 Let $H$ be a spanning subgraph of a graph $G$ which is the
co-support of $\phi$, $i. e.,$ $V(H)=V(G)$ and
$E(H)=\phi^{-1}(id)$, where $id$ is the identity element in $S_n$
and for $E(H)$, we identify each pair of oppositely directed edges
of $\phi^{-1}(id)$. Then our chromatic number of $H$ respect $G$
naturally extends as follows; colorings $f_1, f_2, \ldots, f_n$ of
$H$ are \emph{compatible} if for $e^{+}=(u,v)\in D(G)-D(H)$,
$f_i(u)\neq f_{\phi((u,v)(i)}(v)$ for $i=1$, $2$, $\ldots$, $n$.
The smallest number of colors such that $H$ has an $n$-tuple of
compatible colorings is called the \emph{n-th chromatic number of
$H$ with respect to $G$} and denoted by $\chi_{G}(H)$. Unlike
two fold coverings, estimations of the $n$-th chromatic
numbers of $H$ with respect to $G$ are not easy. The asymptotic behavior of
the chromatic numbers of non-isomorphic $n$-fold coverings could be very fascinating compare to the result by Amit, Linial, and Matousek~\cite{ALM}.

We conclude the
discussion with an example. It is easy to see that all odd-fold
covering graphs of the graph $G$ in Figure \ref{exam1} have
chromatic number $3$. We provide a $4$-fold covering graph induced
by the coloring $\phi$ in Figure \ref{exam3} together with $4$
compatible colorings $f_1$, $f_2$, $f_3$ and $f_4$ of the spanning
subgraph $H= cospt(\phi)$.

\vskip .5cm \noindent{\bf Acknowledgements}

The authors would like to thank Younghee Shin for her attention
to this work. Also, the referees have
been very helpful and critical during refereing and revising. The \TeX\, macro package
PSTricks~\cite{PSTricks} was essential for typesetting the
equations and figures.

\bibliographystyle{amsalpha}

\end{document}